# Learning from Experts: A Survey

Irene Valsecchi
University of Milano-Bicocca (IT)

November 1, 2018

## Abstract

The survey is concerned with the issue of information transmission from experts to non-experts. Two main approaches to the use of experts can be traced. According to the game-theoretic approach expertise is a case of asymmetric information between the expert, who is the better informed agent, and the non-expert, who is either a decision-maker or an evaluator of the expert's performance. According to the Bayesian decision-theoretic approach the expert is the agent who announces his probabilistic opinion, and the non-expert has to incorporate that opinion into his beliefs in a consistent way, despite his poor understanding of the expert's substantive knowledge. The two approaches ground the relationships between experts and non-experts on so different premises that their results are very poorly connected.

Codes: 91A, 91B, 62B

## 1 Introduction

The use of experts is pervasive and problematic. Experts are the common element of very disparate situations, that range from managers who receive technical reports to tourists interested in meteorological forecasts. The problem is that all the agents who are not experts must come to assess the reliability of professional advice in order to use it.

On one side, people need to communicate more and more because knowledge has become very specialized and information widely distributed; on the other side, juts because knowledge is limited and information is segmented, people can be either unable to grasp the opinions of experts fully, or reluctant to trust experts, overall in case of conflicting interests. Consequently, learning from experts may well not be a straightforward matter, and, recently, a steady flow of contributions has come to tackle the problems related to information transmission from experts.

In general terms expertise refers to the case in which one party sends messages to another party who faces uncertainty. The sender of messages is called expert. The receiver of messages is interested in learning from the expert. The decoding of messages and the associated posterior uncertainty of the receiver are the essential issues that go under the heading of information revelation.

The literature about expertise usually resorts to one of the two following approaches: the first approach is grounded on game theory, the latter approach is founded on Bayesian decision theory. In the game-theoretic approach, all the agents, experts and non-experts, share the same *ex ante* view of the world, i.e. they agree upon the space of the states of nature and the space of the realizations that some privately observable variables can take. Overall, they agree upon the correspondence between state and observation in terms of likelihood function of an observation when a particular state of the world is true. The expert usually is the only agent who makes private observations. The game-theoretic approach is mainly concerned with the amount of learning that is compatible with strategic behavior on the expert's side at the communication stage, especially when expert and non-expert have different payoff functions.

Instead, according to the Bayesian decision-theoretic approach, the expert is the sender of his probabilistic opinion about some unknown state. The problems lie in the criteria of learning that the receiver of the expert's opinion should adopt in updating his own beliefs. Here the expert is usually supposed to deliver his probability assessments honestly. His messages are taken to be experimental evidence by the receiver, who has no knowledge of the variables observed by the expert. No strategic interplay between expert and non-expert is assumed: attention is focused on the requirements for the incorporation of the expert's probabilistic opinions.

The contributions generated by those two approaches are rich but poorly exchangeable and overlapping. In the game-theoretic approach the agents need to share the same knowledge of the world so that the expert can play strategies rationalizable by the non-expert. In the decision-theoretic approach expert and non-expert have differential knowledge, but the expert is treated as a non-strategic mechanism generating data. Consequently, the insights of one approach can be looked at as the limits of the other approach. The game-theoretic approach emphasizes that the amount of information transmitted by the expert in equilibrium can be inefficient because the expert can derive advantages from noise in communication. The decision-theoretic approach emphasizes the implications of the non-expert's subjective expectations regarding the reliability of the expert's expertise.

The present survey is focused on problems of information revelation from experts. The concern is on information transmission from experts who formulate probabilistic assessments at no relevant cost. Consequently, the issues of information acquisition on the experts' side are left aside. In models of information acquisition, expertise is effort-expensive, and either contractual[1] or organizational[2] arrangements are studied as a response to incentive problems. Although closely related to information revelation, information acquisition, is omitted because it belongs properly to the literature on principal-agent relationships with unobservable effort and hidden actions.

Also the present survey leaves apart contributions to the economics of credence goods. The key feature of credence goods is that the supplier knows more about the quality of the service a consumer needs not only ex ante but also after the service has been provided. In markets for credence goods the supply side is composed by experts specialized in the diagnosis and the treatment of the consumers' needs. Fraud and cheating, overall in the form of over-treatment, are the major problems in the interaction between experts and consumers. An extensive review about the economics of credence goods is Dulleck and Kerschbamer (2006).

The paper is organized as follows. In Section 2 I survey contributions adherent to the game-theoretic approach. The models are classified according to the role played by the agent who is not an expert. In Section 2.1 the non-expert is a decision-maker who takes some action that is relevant for the expert's payoff. In Section 2.2 the non-expert is an evaluator who measures the expert's performance with payoff-relevant consequences for the expert. Particular attention is devoted to the issue of the expert's concern for reputation. In Section 3 I survey contributions adherent to the Bayesian decision-theoretic approach. In Section 4 I conclude.

Since expertise is connected *per se* to uncertainty and stochastic variables, I will adopt the following notation throughout the paper. The term generalized probability density function ($g.p.d.f.$) will designate a function that may be either a probability function ($p.f.$) or a probability density function ($p.d.f.$). A state of nature is denoted $\omega$, while $\Omega$ is the space of the states of nature, that are mutually exclusive and exhaustive. Agent $i$'s prior $g.p.d.f.$ of $\omega$ is denoted $f_0^i(\omega)$, with distribution function ($d.f.$) denoted $F_0^i(\Omega)$. At time $t$, the $g.p.d.f.$ of $\omega$ is denoted $f_t(\omega)$. An event $A$ is a subset of states of nature. An observation $s$ is a realization drawn from the observation space $\Sigma$. Given that the true state of nature is $\omega$, agent $i$'s conditional $g.p.d.f.$ (or likelihood function) of $s$ is denoted $l^i(s \mid \omega)$.

## 2  Expertise as Private Information

Consider the following situation. An agent called $N$ is uncertain about the true state of nature. An agent called $E$ observes a private and informative signal of the true state. Some communication from agent $E$ to agent $N$ is feasible at no cost, i.e. agent $E$ can send messages to agent $N$. However, messages cannot be verified by agent $N$, possibly because the signal observed by agent $E$ is soft information and no proof can be produced to substantiate agent $E$'s claim. The relevant issue is the precision of the transmitted information, measured by the correspondence between observation and message. The following basic set-up is usually adopted:

*1) Ex-ante uncertainty* (1)

two agents, $N$ and $E$, believe that the true state of nature belongs to the common state space $\Omega$. The agents have prior beliefs about the true state of the world, i.e. each agent has in mind a unique *g.p.d.f.* of $\omega$.

*2) Knowledge* (2)

the agents know that there exists an observation variable $S$ taking values in the observation space $\Sigma$. Both agents agree that a unique complete set $L$ of likelihood functions $l(s \mid \omega)$ (i.e. a single conditional *p.d.f.* of $s$ given $\omega$ for all $\omega$) represents the connection between observations and states.

*3) Information* (3)

agent $N$'s prior *g.p.d.f.* of $\omega$ is non degenerate. Only agent $E$ observes the realization $s$.

*4) Common knowledge* (4)

the agents' prior *g.p.d.f.* of $\omega$, the set $L$ of likelihood functions, the fact that only agent $E$ observes $s$, the agents' payoff functions and the timing of interaction are common knowledge.

*5) Agent E's strategy* (5)

agent $E$'s strategy is a message rule, i.e. a complete set of conditional *g.p.d.f.* of $m$ given $s$.

Agent $E$ is the only agent making an observation directly related to the state of nature. For that reason $E$ is called expert, and expertise is a particular case of asymmetric information, with agent $N$ having strictly less information than agent $E$.

An essential feature of the above set-up is that both agents $E$ and $N$ agree upon the conditional distribution of $s$ given $\omega$ for all the possible $\omega$ (Bayarri-DeGroot, 1991). The supposed consensus upon the connection between observation and state, i.e. a single likelihood function $l(s \mid \omega)$ for each $\omega$, is the necessary common ground for the interplay between Bayesian players.

The relevant economic models can be classified according to the following criterion:

a) agent $N$ is a decision-maker in that he takes some action $a$. The expert's payoff function directly depends on the pair of implemented action and true state.

b) Agent $N$ is an evaluator in that he measures the expert's performance as a predictor.

## 2.1 The expert's concern for the decision-maker's action

Consider the following situation. Agent $N$ is a decision-maker because he can take an action that is payoff-relevant for both himself and agent $E$. However, the same action can yield different payoffs according to the true state of the world. Agent $N$ can listen to agent $E$'s recommendations before deciding which action to implement. The issue is how much information agent $N$ can extract from agent $E$'s messages.

In order to represent the above case, the following specifications usually supplement the set-up in (1)-(5):

*6.a) Choices*

after having observed $s$, agent $E$ chooses one message $m$ out of the message space $M$, and sends it to agent $N$. After having received the message, agent $N$ chooses some action $a$ out of a feasible action space $A$.

*7.a.) Payoff functions*

both the agents' payoff functions depend on the pair $(a, \omega)$.

*8.a) Agent $N$'s strategy*

Agent $N$'s strategy is an action rule, i.e. a complete set of conditional $g.p.d.f.$ of $a$ given $m$.

When agent $N$ is a decision-maker and agent $E$ is privately informed, economic models focus on the welfare properties of the equilibria that characterize the relationship between the agents. In this way, the relationship between decision-maker and expert becomes a particular variant of the class of principal-agent relationships. Crawford and Sobel (1982) are the usual reference paper for the above set-up, in which the expert's payoff depends on the decision-maker's action. In particular, they consider the case in which agent $E$ can perfectly induce the true state of nature, i.e. the observation space $\Sigma$ coincides with the state space $\Omega$, equal to the closed unit interval on the real line. Agent $E$ has a twice continuously differentiable von Neumann-Morgenstern utility function $U^E(a, \omega, b)$, where $a$ is a real number and $b$ is a scalar parameter. Agent $N$ has a twice continuously differentiable von Neumann-Morgenstern utility function $U^N(a, \omega)$. Both payoff functions are strictly concave in action, and the cross derivative $U^i_{as}(\bullet)$ is strictly positive, making each agent's optimal action strictly increasing in $\omega$ under certainty. When the scalar parameter $b$ is not nil, there will be conflict of interest between the agents because the optimal actions of agent $N$ and agent $E$ are always different in each certain state of nature. The focus is on strategic information transmission by the partisan expert, who can be interested in misleading the action of the decision-maker to his own advantage.

In particular, Crawford and Sobel show the following: provided there is conflict of interest, the Nash Bayesian equilibria will imply a partition of the state space into an only finite number of subintervals. In other words, provided the agents' payoffs never coincide for each state of nature, in equilibrium the updated state beliefs of agent $N$, having received message $m$, will have support on a non-singular subinterval of the state space, consistently with agent $E$'s optimal message rule conditional on the realized observation. Moreover, if there is a Nash Bayesian equilibrium implying a partition with $n$ subintervals, there will exist a Nash Bayesian equilibrium implying a partition with less subintervals. Since the Nash Bayesian equilibria are partition equilibria, there will be a finite number of equilibrium actions only. Finally, since the meaning of messages is endogenously determined in equilibrium, any permutation of messages across meanings will yield another equilibrium (Li 2007).

Consequently, there cannot be perfect revelation of private information from agent $E$ to agent $N$ in case of conflict of interest[3]. Some coarsening of information will be unavoidable in equilibrium. Noise in communication is proved to decrease as the preferences of the agents become more aligned. Finally, a Nash Bayesian equilibrium in which no information is disclosed, i.e. a bubbling equilibrium, will always exist. That equilibrium emphasizes the relevance of the non-informed player $N$'s conjectures in bounding the expected welfare from the game: agent $N$'s expectation of the noisiest quality of messages from agent $E$ will always be self-fulfilling, making the expert devoid of any effective role.

Information manipulation, and, overall, babbling equilibria, in which agent $N$ takes all messages as meaningless, are a neat theoretical result. Somehow that result stands in paradoxical relationship with the actual pervasive use of experts. Information disclosure is shown to be enhanced when auxiliary assumptions modify the basic set-up of Crawford and Sobel. In particular, the new assumptions introduce:
  - the delegation of action from the decision-maker to the expert;
  - different modes of communication between expert and decision-maker;
  - incomplete information on the expert's side;
  - iterated interaction between expert and decision-maker;
  - the use of multiple experts instead of a single expert.

*Delegation.* The agency problem can be mitigated by delegating action to the expert. Li and Suen (2004) consider the case in which the action to be taken is the choice of adoption or rejection of a project, that can yield either positive or negative net fixed benefits. As usual, the expert observes signals conditional on the true state. Agents' priors about the true state of nature, agents' preferences and the quality of the observations are common knowledge. Preferences correspond to the losses induced by the wrong decision. Li and Suen show that there will be scope for delegation to the privately informed expert as long as the decision-maker is not *a priori* extremely biased towards a particular choice. In particular, under some conditions, disaligned preferences of agents $N$ and $E$ will actually make delegation (i.e. the commitment not to overrule the expert) credible. The intuition is the following: suppose that agent $N$ is *a priori* biased toward rejection, and that agent $E$ is *a priori* more prone to rejection than agent $N$. If agent $N$ delegates action to agent $E$, agent $E$ will choose according to the his realized observation. If agent $E$ chooses rejection, the expert's choice will correspond to the choice that agent $N$ would have taken given his prior. If agent $E$ chooses adoption, then agent $N$'s updated probability that the project will yield positive benefits will significantly increase, making adoption a welcome decision. In that case, delegation cannot ever be outperformed by agent $N$'s commitment to any other mechanism designed to elicit information from the expert.

*Mode of communication.* The same structure of communication in terms of rounds of messages can affect the amount of information disclosed in equilibrium. Krishna and Morgan (2004) consider the following mode of communication: there is a first round of conversation in which simultaneous and cheap talk messages are exchanged between expert and decision-maker; depending on how the first-round conversation goes, the meeting is deemed by both parties to be either a success or a failure; after this, the expert may send a further message. Random messages from the uninformed decision-maker at the first round of communication can well lead to equilibria that ex ante Pareto dominate the set of equilibria identified by Crawford and Sobel, provided the conflict of interest is not too extreme. The reason is that randomness in conversation breaks the deterministic link between the expert's message and the decision-maker's action, thereby reducing the expert's incentives to withhold information in direct relation with the expert's risk aversion.

*Two-sided incomplete information.* In Crawford and Sobel the expert observes the realization of the payoff relevant state of nature perfectly, while the decision-maker is perfectly uninformed. On the contrary, if there is two-sided incomplete information, the amount of information disclosed in equilibrium will increase. Watson (1996) shows that truthful revelation can hold in equilibrium provided the observations of the expert and the decision-maker are correlated, and it is likely that the decision-maker's information leads to a decision that agent $E$ likes, conditional on his information. When the decision-maker knows the code of the observation privately made by the expert, while the expert does not, then even full revelation can be an equilibrium outcome.

*Repeated interaction.* Instead of static games, dynamic games between agent $N$ and agent $E$ can show how credibility in information transmission is built up. In Sobel (1985) the interaction between agent $N$ and agent $E$ is repeated. Moreover, the decision-maker is supposed to be uncertain about the preferences of the expert, who can be either a "friend" or an "enemy". Iterated interaction between the agents, coupled with verifiable information at the end of each stage of interaction, is shown to make it worthwhile for the sender to build a reputation for truthfulness[4].

*Multiple experts.* The clash between the informed party's partisanship and the revelation of information can be reduced by having multiple experts, who are supposed to have perfect information, or receive noisy signals, or make noisy multidimensional observations in different models. In Krishna and Morgan (2001), two experts observe the true state of nature perfectly, and send public messages to the decision-maker in a sequential fashion. The fully revealing Perfect Bayesian equilibrium is shown never to exist. However, the decision-maker will always derive some benefit from consulting both the experts as long as the bias of one expert's preferences with respect to the decision-maker's ones is opposed to the bias of the other expert's preferences.

Austen-Smith (1993) focuses on imperfectly informed experts. In particular, two biased experts observe conditionally independent signals, given a binary signal space and a binary state space. Different configurations of their preferences relative to one another and the decision-maker are considered. Austen-Smith (1993) shows that sequential reporting is superior to simultaneous reporting in terms of information disclosure, with little difference between like and opposing biases.

Wolinsky (2002) is concerned with the case in which many experts share the same preferences that significantly diverge from the decision-maker's preferences. The decision-maker can choose between two different actions only, while the experts observe signals that are supposed to be identically and independently distributed. Moreover, observations can take only two values, one of which is verifiable. The purpose of the model is to analyze the impact of communication among the experts over the amount of information that the decision-maker can elicit. The mode of communication is the following: the decision-maker is assumed to partition the set of experts into groups; experts can exchange information only if they belong to the same group, while each group can only make a report to the decision-maker. Each group of experts will choose its message with the understanding that it will matter only if pivotal. By allowing communication only within groups, the decision-maker is shown to extract more information than the experts would like to reveal.

Many experts observe noisy signals, conditionally independent given the state of nature, in Gerardi et al. (2005) as well. They show that, when signal are very precise and messages are aggregated, an elicitation scheme implying randomized actions can make no single expert expect to change the implemented action by misreporting his information. For a similar effect, information can be extracted at small cost in terms of deviation from the optimal action rule, even if each expert's information is not very accurate but the number of consulted experts gets large.

Battaglini (2002) proves that the dimensionality of the state space and the action space has an important impact on the transmission of information in equilibrium with many experts. He shows that full revelation is possible when there are at least two experts who are perfectly informed, provided the state space is multidimensional, even with an arbitrarily large conflict of interest. On the contrary, in Battaglini (2004) noise in the observations made by many experts leads to the impossibility of full revelation of information, independently of the number of experts or the dimensionality of the state space. Experts are supposed to observe unbiased signals, while noises in the observations made by different experts are independent. The reason of the impossibility of full information revelation is that noise generates a trade-off between information aggregation and information extraction. Since signals are noisy, it is important to aggregate observations. However, the conflict of interest between experts and decision-maker will lead the decision-maker to limit the dimensionality of the messages.

To sum-up, the amount of information that a decision-maker can extract is shown to increase somehow through the aggregation of reports in a sequential or in a parallel mode. In particular, the use of multiple experts calls for strict assumptions regarding the joint density function of the observations of the different experts. No attention is paid to the decision-maker's costs from consulting many experts, for instance in the form of optimal stopping rules in the process of search, information overload or delay in reaching a final decision.

## 2.2 The expert's concern for evaluation

The expert's concern for evaluation plays a prominent role when agent $N$ measures the expert's expected or current performance as a predictor. Most of the literature focuses on measures of expected performance, i.e. on the expert's concern for reputation.

Consider the following situation. All agents are uncertain about some relevant parameter of the environment. Agent $E$ comes from a population of privately informed agents, endowed with different levels of ability as non-strategic forecasters of the environmental parameter. That is, the informative property of the signal, that is observed by agent $E$ and that concerns the environmental parameter, depends on some personal characteristic of the same agent $E$. *Ex ante* agent $N$ does not know the ability of agent $E$. Agent $N$ has just prior beliefs about the distribution of ability in the population of the privately informed agents. Agent $E$ sends one message to agent $N$, who comes to observe the true environmental parameter after the message has being received. Agent $E$ is concerned with his final reputation for being an high-quality predictor in the opinion of agent $N$, and not with the decisions that can be implemented because of his advice. In this sense agent $N$ is an evaluator of the quality of agent $E$ as a forecaster.

In order to represent the above case, the following specifications usually supplement the basic set-up in (1)-(5):

*1.1) Ex-ante uncertainty*

the state of nature $\omega$ is the pair $(x, \eta)$. The realization $x$ represents some parameter of the environment, while the realization $\eta$ represents some characteristic of agent $E$. Both agents, $E$ and $N$, believe that the realization $x$ belongs to the common parameter space $X$, while the realization $\eta$ belongs to the common space $\Xi$ of characteristics. The state space $\Omega$ is the Cartesian product of the spaces $X$ and $\Xi$.

*2.1) Knowledge*

Both agents believe that $x$ and $\eta$ are statistically independent. The prior $g.p.d.f.$ of $\omega$ of each agent $i$ is the product of the $g.p.d.f.$ $g_0^i(x)$ and the $g.p.d.f.$ $p_0^i(\eta)$.

*3.1) Information*

At the end of the period the true parameter $x$ becomes public information.

*6.b) Choices*

after having observed $s$ but not $x$, agent $E$ chooses one message $m$ out of the message space $M$, and sends it to agent $N$. At the end of the period, after having observed both $m$ and $x$, agent $N$ computes his posterior belief $p_1^N(\eta \mid m, x)$.

*7.b) Payoff functions*

agent $E$'s payoff function depends on $p_1^N(\eta \mid m, x)$.

*8.b) Agent N's strategy*

agent $N$'s strategy is a Bayesian updating rule.

The paper of Ottaviani and Sorensen (2006a) is representative of the above set-up, suited to the analysis of the relationship between truthful reporting of information and reputational concern. They assume that both the space of environmental parameters $X$ and the space of personal characteristics $\Xi$ are intervals on the real line. All the agents share common and non-degenerate prior beliefs on the true state so that they are uncertain about both the realization $x$ and the realization $\eta$. The personal characteristic $\eta$ can be interpreted as ability in forecasting because agent $E$'s observation is generated by an experiment for which the conditional $p.d.f.$ $l(s \mid x, \eta)$ in the set $L$ can be ordered according to $\eta$. That is, $l(s \mid x, \eta)$ is a sufficient experiment for $l(s \mid x, \eta')$ when $\eta$ is higher than $\eta'$: a higher realization $\eta$ corresponds to better information in the sense of Blackwell. In this way, the personal characteristic $\eta$ parametrizes the amount of information about $x$ contained in agent $E$'s signal. In particular, agent $E$'s observation is supposed to be a multiplicative linear experiment, i.e. a mixture of an informative and an uninformative component. A higher $\eta$ is more likely to make the observed signal drawn from the informative component. The likelihood function $l(s \mid x, \eta)$ is assumed to satisfy the monotone likelihood ratio property[5] in $(s, x)$ for any given $\eta$. Agent $E$'s preferences are represented by his reputational payoff from message $m$ when the true environmental parameter is $x$. That reputational payoff corresponds to $\int_\Xi v(\eta) p_1(\eta \mid m, x) d\eta$, where $v(\eta)$ is a strictly increasing function. The justification for agent $E$'s reputational payoff is that reputation for ability in forecasting is rewarded by non-modelled decision-makers. Agent $E$'s expected payoff from sending message $m$ is shown to be strictly supermodular[6] in the pair $(m, s)$. In this sense, for agent $E$ each message $m$ will induce a lottery over posterior reputations with final payoffs that depend on the true environmental parameter. Agent $E$ will choose which message to send according to his posterior beliefs about the environmental parameter, updated given his observed signal.

Ottaviani and Sorensen (2006a) show that agent $E$ will gain from false reporting if agent $N$ conjectures truthful reporting. Consequently, the message maximizing agent $E$'s expected reputational payoff is typically different from the signal actually observed[7]. As in the case of an expert concerned with a decision-maker's action, weak Bayesian Nash equilibria, that impose no restrictions on agent $N$'s beliefs for off-path messages, are proved to have a partition structure: each message sent in equilibrium corresponds to signals in some subset of the observation space. Once again, there will always exist a babbling equilibrium, but now that equilibrium will occur when the prior belief on the environmental parameter is sufficiently concentrated, and not when the conflict of interest is sufficiently severe. Ottaviani and Sorensen (2006a) show that the most informative equilibrium is binary, i.e. at most two messages are transmitted in equilibrium. It follows that agent $E$ can at best communicates the direction but not the intensity of his information[8].

Similarly to the case of the expert's concern for the decision-maker's action, so the impact of the expert's reputational concern on efficiency is analyzed under different specifications of the basic set-up. In particular, the following modifications are inquired:

- many experts exchange either simultaneous or sequential messages;
- the temporal horizon of the expert's activity is extended.

*Many experts.* Ottaviani and Sorensen (2006b) analyze whether the report-ing incentives of experts improve once reputational concern is coupled with competition among many forecasters. The environmental parameter $x$ and personal ability $\eta$ are independently distributed, as well as the experts' signals conditional on $x$ and $\eta$. Experts care about their relative reputation: each expert $i$ is assumed to have a von Neumann-Morgenstern payoff, that depends on both the ability of agent $i$ and the ability profile of all the other experts. Ottaviani and Sorensen prove that in equilibrium, with experts either not knowing or knowing their own ability, each expert will behave as in the absolute reputation model. The reason is that the independence of the experts' signals implies stochastic independence of the posterior reputation of different experts, updated after the reports and the true environmental parameter become public information. Hence, only an expert's own message can influence the updating of his reputation.

Scharfstein and Stein (1990) consider the case of multiple experts who receive correlated signals and send sequential messages. They assume that both the parameter of the environment and the personal characteristic of an expert are binary. Two experts send public messages in a sequential mode. In particular, an expert's message is taken to be a public and dichotomous investment decision, that will be profitable only if the high environmental parameter is true. An expert's observation is assumed to be imperfectly informative of the true environmental parameter only if the expert's ability is high, otherwise the signal is perfectly uninformative. *Ex ante* the distribution of signals is the same for high and low ability experts. If the two experts have either low ability or different abilities, their observations will be stochastically independent; instead, if both the experts have high ability, their observations will be perfectly correlated, i.e. they will observe the same signal. Agent $N$ is an evaluator in that, after both the experts have sent their messages, he will update beliefs about the experts' ability based on two pieces of evidence, a) the true environmental parameter and b) both the experts' messages. No expert knows his own ability. Each expert is supposed to care about his posterior ability reputation. Scharfstein and Stein show that herd behavior can hold in equilibrium because of a "sharing-the-blame" effect. In particular, the expert who sends his message after having observed the message of the previous expert can mimic the other expert's message, ignoring his private substantive information, just because all high ability agents receive identical signals[9].

Ottaviani and Sorensen (2001) analyze herding in committees of experts when the environmental parameter, ability, observation and message are all binary. The distinguishing feature is that, conditional on the environmental parameter, the observations made by different experts are supposed to be statistically independent. Again all the agents have common prior beliefs on ability, and a high-ability expert receives a more precise, although imperfect, signal. Given a group of experts commonly known to have heterogenous abilities, each agent $E$ is supposed to send a public message only once according to a pre-ordered sequence of messages. After all the messages have been sent and the true environmental parameter has become public information, agent $N$ updates the reputation of each expert. Each expert's payoff is assumed to be linearly increasing in his posterior reputation so that an expert's objective is to maximize the expected value of his reputation, conditional on his private information. Ottaviani and Sorensen prove that experts may herd by suppressing their true observations, although stochastically independent. The extraction of information can be improved only up to some imperfect degree by reordering the sequence of messages.

Instead, the exchange of private information is simultaneous among the experts in the committees analyzed by Visser and Swank (2007). Again the environmental parameter, ability and observation are binary. Each expert does not know his own ability; however, if he has high ability, his observation will be perfectly informative of the environmental parameter, while, if he has low ability, his observation will be perfectly uninformative. The committee of experts is supposed to reach a decision in two stages: in the first stage, each expert can exchange messages with the other experts; in the second round, the experts vote simultaneously and a final public decision is reached on the basis of some majority voting rule. Given the final decision, agent $N$ updates his beliefs about the ability of each expert. Expert $i$'s payoff is assumed to depend on both the revenues from the final decision and agent $N$'s posterior probability that expert $i$ has high ability. While experts share common preferences as to the final decision, they differ in the weight they put on their reputational concern. Visser and Swank prove that reputational concerns will lead the committee to show a united front and may distort the final decision. If the experts are very heterogeneous in their preferences for reputation, they will manipulate information at the communication stage[10].

*Repeated forecasting.* Prendergast and Stole (1996) are concerned with the impact of reputational concern on the use that the expert makes of his own private information over time. In each period $t$ agent $E$ sends a public message, after having made an observation in the same period. Agent $E$'s message at time $t$ is interpreted as the amount of money invested in a project at time $t$. The project revenues at time $t$ are linear in one unknown environmental parameter $x$, constant in all periods, and concave in agent $E$'s current investment. The prior and common distribution of the environmental parameter $x$ is normal with zero mean and known variance. At the beginning of each period, agent $E$ makes a private and imperfect observation of the true environmental parameter; in particular agent $E$'s observation at time $t$ is supposed to be distributed as the sum of the true environmental parameter and a current error component. The current error component has zero mean and variance that can take only two fixed values according to agent $E$'s ability. Ability is a personal characteristic known to agent $E$, but *ex ante* it is drawn from a population with either high or low ability in making observations. Consequently a high-ability expert will receive more precise estimates of the true environmental parameter $x$. Agent $E$'s expected payoff in period $t$ depends on both his posterior expected current project revenues and his immediate end-of-period reputation for ability. End-of-period reputation is measured as agent $N$'s equilibrium expectation of agent $E$'s observational variance, updated on the profile of agent $E$'s messages from period 1 up to period $t$. Since inferences are drawn from the difference between actual and expected investment, the variance of the posterior estimate of the environmental parameter will play an important role in agent $N$'s inference of agent $E$'s ability. Prendergast and Stole (1996) show that initially agent $E$ will overreact to his private new information, consistently with an overconfidence effect. Indeed, at the beginning, efficient investment is more variable for high-ability experts. Since the variance of the expert's posterior is increasing in ability, agent $E$ will exaggerate his true beliefs about $x$ on the margin in order to appear more talented. However, after some time, agent $E$ will become unwilling to respond to his private new information because otherwise agent $N$ would induce that his previous behavior was wrong. That later behavior is compatible with conservatorism or sunk-cost fallacy. The reason is that later messages are an indication of the precision of agent $E$'s both previous and current observations, and a high-ability expert should change little. Hence, provided informational gains from additional observations are sufficiently decreasing, there will follow a too little variability in investment.

Sequential forecasting with agent $N$ being both a decision-maker and an evaluator is the concern of Li (2007). In a two-period model, the environmental parameter, ability and action are binary. Environmental parameters are independent across time. Agent $E$ knows his ability and makes two observations in each period, before the environmental parameter becomes public information at the end of each period. Agent $E$'s observations are independent conditional on the environmental parameter, with the second observation being more accurate than the first observation. Moreover, all the observations of a high-ability expert are more precise than the observations of a low-ability expert. Agent $E$ is assumed to be motivated only by his second period wage that is increasing in agent $N$'s posterior belief about his ability. So agent $E$ is driven only by reputational concerns, while agent $N$ cares about the revenues yielded by the implemented action and the payments made to agent $E$. For the case in which agent $E$ is required to send a message after each observation, Li shows that in perfect Bayesian equilibria inconsistent reports within the same period can signal ability, provided that the second observation of a high-ability expert improves faster than the second observation of a low-ability expert. Contrarily to Scharfstein-Stein (1990) and Prendergast-Stole (1996), consistent reports *per se* may indicate low ability in equilibrium. Finally, the sequential reporting protocol will be more information revealing than a single message if the agent $N$'s action is very sensitive to the accuracy of the reports.

As well as the expert's concern for his expected performance, the expert's concern for his current performance as a predictor can lead to strategic forecasting. Ottaviani and Sorensen (2006c) analyze a forecasting contest with pre-specified rules. All the agents share the common prior belief that the environmental parameter is normally distributed with known mean and precision. An expert's private observation is normally distributed around the true environmental parameter. Conditional on parameter $x$, the experts' signals are independently normally distributed with mean $x$ and common precision. After having made their observation, all the experts simultaneously submit their forecasts, i.e. their *ex-post* expected environmental parameter $x$. Once the true environmental parameter becomes public information, agent $N$ acts as an evaluator and uses the experts' forecasts and the realized environmental parameter to evaluate forecasting ability. In particular, it is assumed that the expert whose report is closest to the true $x$ wins a prize to be divided among all the winning experts. For the case of an infinite number of experts, Ottaviani and Sorensen prove that in the symmetric Bayes-Nash equilibrium experts will put greater weight on their private signals than they would in a honest report of their posterior expectations. The reason is that two opposing forces drive the experts' behavior. On one hand, an expert has an incentive to report his honest forecast, which is most likely to be correct. On the other hand, a forecaster gains from moving away from the prior mean because the number of experts correctly guessing the parameter will be lower if the environmental parameter is away from the prior mean.

All the above contributions make the recognition of different forecasting capabilities their starting point. Not only experts can introduce noise into their messages strategically, but their announcements can be affected by different biases because they originate from observations of different informative content. However, although justified by the connection between the experts' rewards and use in the background, reputational concerns remain only a proxy of the actual role played by the experts' recommendations in taking decisions. Again the assumptions regarding the joint density function of either the observations of the single expert over time or the observations of many experts have a strong impact on the results of the models.

## 3  Expertise as probabilistic opinions

The Bayesian approach to the use of an expert's advice implicitly posits that the hiatus between an expert and a non-expert stretches beyond some gap in private information: the decision-maker does not share the same knowledge of the expert. In particular, the decision-maker cannot deconstruct the expert's probability assessment into the evidence supporting it. No correspondence between states of nature and observations of the expert is postulated. Hence the Bayesian approach leaves apart the issue of strategic information transmission because it leaves apart the assumption according to which there exists an objective $g.p.d.f.$ of the expert's observations conditional on the state of nature. Consider the following set-up:

*1) Ex-ante uncertainty*

two agents, $N$ and $E$, believe that the space of the states of nature can be partitioned into the event $A$ and the event non-$A$, denoted $\overline{A}$. The agents have initial subjective probabilities about those events, respectively $f_0^i(A)$ and $f_0^i(\overline{A})$.

*2) Information*

agent $N$'s prior probability of event $A$ is non-degenerate. Agent $E$ announces his subjective probability $f_0^E(A)$ to agent $N$ truthfully.

*3) Choice*

having observed $f_0^E(A)$, agent $N$ must compute his posterior belief $f_1^N\left(A \mid f_0^E(A)\right)$.

The relevant problem is the following: how should the decision-maker alter his prior upon reception of $f_0^E(A)$? The problem is not trivial once we consider that the background of the expert problem is a situation of significant uncertainty and complexity according to French (1986). The decision-maker has little substantive knowledge of the factors affecting the event of interest, thereby asking another person for advice. The consulted person is referred to as an expert. In this contest an expert is anyone who can give predictions, i.e. anyone who can make probability statements, called judgments or opinions, concerning the event of interest. In particular, in case of events essentially unique, the assessment of their probability can only be based on personal judgements, as Lindley et al. (1979) point put. Hence the question is: what rules should the decision-maker use in order incorporate an expert's opinion into his own?

According to Morris (1974) consulting an expert is like performing an experiment. Just as the results of an experiment are *a priori* unknown to an experimenter, so the expert's advice will be uncertain to the decision maker prior to receiving it. The decision-maker should look upon the expert's opinion simply as a piece of data, i.e. the expert's opinion should be treated as a random variable the value of which is revealed to the decision-maker (Genest-Schervish, 1985). Indeed, Morris (1977) underlies that a distinction is required between the meaning of an expert's probability assessment to the decision maker and to the expert himself: the expert views his probability assessment as a reflection of his state of information, while the decision-maker looks at the expert's probability assessment as information itself. In this way, the expert use problem becomes one of Bayesian inference. The only restriction to the posterior assessment of the event $A$ by the decision-maker, conditioned on receiving the expert's opinion, is that it is consistent with both his prior knowledge about the variable and his appraisal of the expert. Specifically, the decision maker must make a subjective appraisal of the dependence between the expert's advice and the actual event. That is, the decision maker must ask himself what his assessment of $f_0^E(A)$ would be if an honest clairvoyant told him that event $A$ will occur. Applying Bayes' theorem, for any given value of $f_0^E(A)$ the decision-maker can calculate his posterior probability of event $A$ to be:

$$f_1^N\left(A \mid f_0^E(A)\right) = \frac{l\left(f_0^E(A) \mid A\right) f_0^N(A)}{l\left(f_0^E(A) \mid A\right) f_0^N(A) + l\left(f_0^E(A) \mid \overline{A}\right) f_0^N(\overline{A})}$$

In particular, $l\left(f_0^E(A) \mid A\right)$ is a likelihood function, subjectively specified, to be interpreted as the probability of the event that the expert's prior is $f_0^E(A)$, given that event $A$ is true. That likelihood function is the model of the expert in the decision-maker's view, specifically it is the decision-maker's subjective measure for the credibility of the expert. If the decision-maker believes the expert to be highly competent, the likelihood function will be strongly dependent on $A$, and it will alter the decision-maker's prior belief significantly. If the decision-maker believes the expert to be highly incompetent, the decision-maker's updated belief will tend to be equal to his prior belief, so invariant to the expert's advice.

Hence, in order to update his beliefs, the decision-maker must proceed to elicit his beliefs about the expert's experience and information thoroughly. This requirement brings about quite problematic issues that affect the Bayesian approach to the use of experts. First, the choice of a model can be complicated by the lack of empirical evidence concerning the expert's performance as a predictor, as well as by a desire for tractable solutions. For that reason, Genest and Schervish (1985) consider the case in which initially the decision-maker is willing to specify only certain features of the joint distribution of the event of interest $A$ and $f^E(A)$. They show how a formula can be chosen which is a coherent posterior probability, no matter what is the decision-maker's complete marginal distribution of $f^E(A)$. In particular, if and only if the agent $N$'s posterior belief is a linear function of $f^E(A)$, it will satisfy a consistency condition, requiring that there exists a joint probability distribution of event $A$ and $f^E(A)$ compatible with that posterior belief[11]. The issue is further developed by West-Crosse (1992) and West (1992). Additional difficulties arise when uncertainty concern non-dichotomous events. In that case the expert is usually assumed to supply a complete subjective probability measure for the uncertain variable of interest. Gelfand et al (1995) are concerned with suitable probability densities for data that are the expert's opinion in the form of an only partial probabilistic specification.

Second, the expert's expertise as a probability assessor is related to the idea of calibration. The decision-maker must consider potential overstatements or understatements of probabilities by the expert (West-Crosse, 1992). Frequency calibration pertains to the agreement between an expert's predictions and the actual observed frequency of the event of interest (De Groot-Fienberg, 1983). If an expert assigns probabilities to events sequentially, an expert is frequency-calibrated if the long run frequency of event $A$ is $f^E(A)$ among those times for which his prediction is $f^E(A)$. Since the decision-maker is an outside assessor, Lindley (1982) proposes the concept of probability calibration: an expert will be probability calibrated if the decision maker adopts the expert's opinion for his own.

A further problematic issue is the relation between the expert's advice and the decision-maker's information, in terms of conflicting or common information. Indeed, if the expert were the ultimate yes-man, the decision-maker might modify his opinions upon hearing them replayed (French 1980, French 1986). The relevance of the problem of information dependence among probability pre-dictors is amplified when many experts supply their forecasts (Clemen-Winker, 1990). A related and relevant issue is that of consensus: given the subjective probability assessments of many experts, how can a single common distribution for the variable of interest be reached? Consensus can be looked at as a problem of aggregation of many opinions (Winkler, 1968, Clemen-Winker, 1999), or dialogue among experts (DeGroot, 1974, Bacharach, 1979). Genest and Zidek (1986) survey the problem of combining probability distributions extensively.

Finally, there might be some difficulty with the assumption according to which the experts deliver their opinion honestly. Bayarri and DeGroot (1989) consider the case in which a finite group of experts quantify their opinions in the form of probability distributions of the state of nature, the true value of which will be public information in the future. Given a finite state space, agent $N$ combines the reported probability functions according to some rule that assigns non-negative prior weights to the opinions of the different experts. After the true $\omega$ has been observed publicly, the weights are updated according to how well the true $\omega$ has been predicted by each expert. Since the updated weights are going to be used in the future, the expert can deviate from reporting his honest subjective distribution in order to maximize his expected weight. Expert $i$'s optimal report will put probability one on a particular state of nature when that his prior weight is close to zero. Instead, when the expert's prior weight is close to one, he will never make dramatic reports. Bayarri and DeGroot assume that each expert $i$ can express his belief about the other experts' report. In particular, each expert $i$ has a likelihood function of the report of all the experts other than $i$, conditional on the state of nature. Unlike the game-theoretic approach, that likelihood function refers to the reports, and not to the honest observations, of the other experts. When expert $i$ is uncertain about the others' report, he may expect to lose weight if he believes that the other experts are good forecasters, who are likely to assign a high probability to the actual state of the world. Hence expert $i$ can report probability zero for a state of nature, even though his honest subjective probability for that state is positive. Moreover, the only differentiable utility function of the expert for which honesty is always the best policy is proved to be the logarithm of the ratio between expert $i$'s posterior weight and the others' posterior weights.

According to the Bayesian approach, the decision-maker looks at the expert as an experiment, as a data generating mechanism. In this sense, the distinction between risk and ambiguity does matter. Epstein and Schneider (2007) generalize the Bayesian learning model in order to accommodate changes in both beliefs and confidence in probability assessments. Their contribution may be relevant for the use of an expert's opinion, although the authors do not consider expertise directly. A decision-maker holds the *a priori* view that his observations are generated by the same memoryless mechanism every period. He deals with a triple of sets, $\Theta$, $P_0$ and $L$. As in the usual Bayesian approach, $\Theta$ is the parameter space, and it represents the features of the data generating mechanism that the decision-maker views as constant over time and that he expects to learn. Instead, $P_0$ is a set of priors (i.e. probability measures on $\Theta$), and it represents the initial view of the decision-maker about the parameters.

When $P_0$ is not a singleton, the decision-maker lacks confidence in the prior information upon which his initial beliefs are based. Finally, in every period the decision-maker makes an observation drawn from a finite observation space. $L$ is a set of likelihoods, and it represents the agent's *a priori* view of the connection between observation and true parameter. A set of likelihoods, instead of a single likelihood for each pair of observation and parameter, is a way of representing the decision-maker's awareness of how poor his understanding is of other factors affecting his observations. Conditional independence is assumed so that past observations affect beliefs about future observations only to the extent that they affect beliefs about the state. A theory about how the observation sample was generated is a pair $(p_0, l^t)$ where $p_0$ is a prior belief on $\Theta$ and $l^t = (l_1, , , l_t)$ $\in L$ is a sequence of likelihoods. How well a theory explains the observations $s^t = (s_1, ..., s_t)$ is captured by the unconditional observations density evaluated at $s^t$, i.e. $\int \prod_{j=1}^{t} l_j(s_j \mid \theta) \, dp_0(\theta)$. A posterior probability measure on $\Theta$ can be derived from each theory by Bayes rule given $s^t$. Reevaluation will take the form of a likelihood ratio test: the decision-maker will discard all theories that do not pass a likelihood ratio test against an alternative theory that put maximum likelihood on the sample[12].

Ambiguity is the essential feature of the model of Olszewski and Sandroni (2007). A decision-maker who does not know the payoff-relevant probabilities for his action can offer a contract to an expert in order to receive his probabilistic opinion. The contract has that the experts communicates the relevant probabilities immediately, while wealth transfers are conditional on the true state, once it has become public information at the end of the period. The expert can be either informed or uninformed. If the expert is informed, he will know the relevant probabilities almost perfectly; otherwise, he will know nothing. State and action are binary, while the informed expert is supposed to announce his probabilistic distribution truthfully. Olszewski and Sandroni adopt the maxmin expected utility representation, according to which the value of an uncertain prospect is determined by the minimum of the expected utilities over a class of possible probability measures à la Gilboa-Scheimedler (1989). They show that there does not exist a screening contract that can separate the informed experts from the uninformed ones.

The Bayesian decision-theoretic approach focuses on the decision-maker's beliefs and subjective measures of the expert's reliability. The same problematic issues that affect the decision-theoretic approach spring from renouncing to postulate an explicit relationship between state and observation of the expert. The introduction of ambiguity reduces the distance between the game-theoretic approach and the decision-theoretic approach. Indeed, the non-expert's a priori view of the connection between observation and state, represented by a set of likelihoods, can be looked at as an alternative specification of different types of experts motivated by reputational concerns. Moreover, the use of one-step ahead conditionals makes the analysis of sequential learning more tractable. Still the gap between the two approaches is unabridged.

# 4 Conclusions

Game theory and Bayesian decision theory are the methodological tools applied to the problem of the use of experts.

Bayesian equilibria emphasize the crucial role played by the non-expert's conjectures in bounding the expected welfare from the expert's recommendations. Indeed, babbling equilibria are the common result of most game-theoretic models. However, in the game-theoretic contributions the appraisal of the expert by the non-expert is disciplined by the assumption of common likelihood functions, that represent the connection between state of nature and observation of the expert.

Instead, according to the Bayesian decision-theoretic approach, the non-expert needs to measure the expert's reliability subjectively. In that way, the non-expert computes his subjective model of the expert, and in that task the non-expert may well be unaided by any substantive knowledge of the expert's observations. Opinions[13] play a central role because the expert announces probability assessments and the non-expert elicits his beliefs about the expert's credibility. However, the strategic element of the interaction between the expert and the non-expert is left apart: the expert is treated as a non-strategic data generating mechanism.

The results achieved by both the approaches are rich. Yet the reconciliation of the substantial divide between them waits for new contributions able to couple subjective learning behavior of non-experts with strategic behavior of experts.

# 5 Notes

(1) Information acquisition from the experts and contractual design is the subject of many papers stemming from Demski-Sappington (1987), Osband (1989) and Prendergast (1993). Information acquisition and management style are analyzed by Rotemberg-Saloner (1993). While in the above literature the expert is modelled as the agent in a principal-agent relationship, leadership is the case in which an informed principal needs to induce agents to exert effort by persuasion (e.g. Hermalin, 1998).

(2) Information acquisition and delegated action is the subject of a growing literature. Aghion-Tirole (1997) and Stein (2002) are reference papers for organizational design as an incentive to exert effort in information acquisition.

(3) In Crawford-Sobel perfect information transmission will be an equilibrium outcome if the agents share the same preferences. Caplin and Leahy (2004) are concerned with an agent $E$ empathetic with agent $N$. Information transmission is still troublesome because agent $N$ can have different preferences over the timing of resolution of uncertainty.

(4) Honesty is analyzed by Olszewski (2004) in one-shot cheap-talk games. The decision-maker is supposed to make an observation non independent from the expert's observation, and the expert's payoff function is assumed to incorporate a concern for reputation.

(5) The likelihood ratio $[l(s \mid x, \eta) / l(s \mid x', \eta)]$ is increasing in $s$ for $x$ greater than $x\prime$.

(6) $\left[\int_X \int_\Xi v(\eta) p_1(\eta \mid m, x) g_1(x \mid s) d\eta dx - \int_X \int_\Xi v(\eta) p_1(\eta \mid m', x) g_1(x \mid s) d\eta dx\right]$ is strictly increasing in $s$ when $m$ is greater than $m\prime$.

(7) Strategic forecasting, induced by measures of the expert's performance based on posteriors, is analogous to strategic signaling from firms concerned with their expected share price, instead of their expected profits. Branderburger and Polak (1996) consider the case in which environmental variable, observation and message are all binary. The expert is a firm that will realize high profits if the high (low) message is sent when the true environmental parameter is high (low). A firm is supposed to care about its expected share price. The equilibrium price of the firm is equal to the posterior probability that the firm has made the correct decision from the point of view of an evaluator. The share price maximization is proved to lead to relevant information being wasted in equilibrium.

(8) For the case in which agent $E$ knows his own type, Ottaviani and Sorensen (2006b) show that a binary equilibrium exists for any prior belief about the environmental parameter $x$.

(9) Under very similar assumptions about the parameter and ability spaces, as well as about the distribution of the experts' observations, strategic forecasting of multiple experts is studied by Laux and Probst (2004). They are concerned with recommendations from experts to decision-makers, and with experts' preferences depending on the ranking of their past performance.

(10) In the above contributions the set of the experts is well identified and common knowledge. Instead, Piccione and Tan (1996) analyze the equilibrium bidding behavior in first-price sealed-bid auctions when bidders are uncertain about how many rivals are experts (i.e. have acquired private information) for both private and common value auctions.

(11) Specifically, $f_1^N(A \mid f_0^E(A)) = f_0^N(A) + \lambda(A)\left[f_0^E(A) - Z\right]$ where $Z = E^N\left[f_0^E(A)\right]$. Nakata (2003) interprets the above result in the following way: when agent $N$ merely computes $f_0^N(A)$ and $Z$, then his posterior belief will be a linear function of $f_0^E(A)$ if and only if it satisfies a razionalizability condition.

(12) Epstein and Schneider adopt a model of recursive multiple priors utility, where beliefs and confidence are determined by a set of one-step-ahead conditionals. The issue of incomplete knowledge about the quality of observations is extended by Epstein and Schneider (2008).

(13) In the economic literature, the issue of opinions, and, more specifically, of differences of opinions is often related to the discussion of the common prior assumption and to differential interpretations of public signals. For instance, in Harris and Raviv (1993), traders update their beliefs about an asset's returns using their own likelihood function of the relationship between public signals and the asset's returns. The fact that traders adopt different likelihood functions is supposed to be common knowledge.